\documentclass[twocolum,fleqn]{autart}    

\usepackage{graphicx}          

\usepackage{amssymb}
\usepackage[cmex10]{amsmath}
\usepackage{accents}
\usepackage{epsfig}

\usepackage{bm}
\usepackage{amsmath}
\usepackage{mathptmx}
\usepackage{color}
\usepackage{dsfont}
\usepackage{empheq}
\usepackage{hyperref}
\usepackage{array}
\usepackage{cite}
\allowdisplaybreaks
\usepackage{subfig}
\usepackage{pgfplots}
\usepackage{tikz}
\usetikzlibrary{arrows,shapes,patterns,positioning,calc}
\tikzstyle{boxstyle}=[draw=black,inner sep=7pt]     
\tikzstyle{arrowstyle}=[]    
\usetikzlibrary{arrows, decorations.markings}

\newlength\figureheight
\newlength\figurewidth

\newtheorem{coro}{Corollary}
\newtheorem{remark}{Remark}
\newcommand{\pa}{\partial}
\newcommand{\ba}{\begin{align}}
\newcommand{\ea}{\end{align}}
\newcommand{\fr}{\frac}
\newcommand{\gm}{\gamma}
\newcommand{\ep}{\varepsilon}

\begin{document}

\begin{frontmatter}

\title{Input-to-State Stability for the Control of Stefan Problem with Respect to Heat Loss} 
 

\author[UCSD]{Shumon Koga}\ead{skoga@eng.ucsd.edu},               
\author[iasson]{Iasson Karafyllis}\ead{iasonkar@central.ntua.gr}, 
\author[UCSD]{Miroslav Krstic}\ead{krstic@ucsd.edu}  

\address[UCSD]{Department of Mechanical and Aerospace Engineering, University of California, San Diego, La Jolla, CA 92093-0411 USA}                                               
\address[iasson]{Department of Mathematics, National Technical University of Athens, 15780 Athens, Greece}
          

\begin{abstract}                          
This paper develops an input-to-state stability (ISS) analysis of the Stefan problem with respect to an unknown heat loss. The Stefan problem represents a liquid-solid phase change phenomenon which describes the time evolution of a material's temperature profile and the liquid-solid interface position. First, we introduce the one-phase Stefan problem with a heat loss at the interface by modeling the dynamics of the liquid temperature and the interface position. We focus on the closed-loop system under the control law proposed in \cite{Shumon16} that is designed to stabilize the interface position at a desired position for the one-phase Stefan problem without the heat loss. The problem is modeled by a 1-D diffusion Partial Differential Equation (PDE) defined on a time-varying spatial domain described by an ordinary differential equation (ODE) with a time-varying disturbance. The well-posedness and some positivity conditions of the closed-loop system are proved based on an open-loop analysis. The closed-loop system with the designed control law satisfies an estimate of ${L}_2$ norm in a sense of ISS with respect to the unknown heat loss. The similar manner is employed to the two-phase Stefan problem with the heat loss at the boundary of the solid phase under the control law proposed in \cite{koga2018CDC}, from which we deduce an analogous result for ISS analysis. 
\end{abstract}

\end{frontmatter}

\section{Introduction}

Liquid-solid phase transitions are physical phenomena which appear in various kinds of science and engineering processes. Representative applications include sea-ice melting and freezing~\cite{Shumon17seaice}, continuous casting of steel \cite{petrus2012}, cancer treatment by cryosurgeries \cite{Rabin1998}, additive manufacturing for materials of both polymer \cite{koga2018polymer} and metal \cite{chung2004}, crystal growth~\cite{conrad_90}, lithium-ion batteries \cite{koga2017battery}, and thermal energy storage systems \cite{zalba03}. Physically, these processes are described by a temperature profile along a liquid-solid material, where the dynamics of the liquid-solid interface is influenced by the heat flux induced by melting or solidification. A mathematical model of such a physical process is called the Stefan problem\cite{Gupta03}, which is formulated by a diffusion PDE defined on a time-varying spatial domain. The domain's length dynamics is described by an ODE dependent on the Neumann boundary value of the PDE state. Apart from the thermodynamical model, the Stefan problem has been employed to model several chemical, electrical, social, and financial dynamics such as tumor growth process \cite{Friedman1999}, domain walls in ferroelectric thin films \cite{mcgilly2015}, spreading of invasive species in ecology \cite{Du2010speading}, information diffusion on social networks \cite{Lei2013}, and optimal exercise boundary of the American put option on a zero dividend asset \cite{Chen2008}.

 While the numerical analysis of  the one-phase Stefan problem is broadly covered in the literature, their control related problems have been addressed relatively fewer. In addition to it, most of the proposed control approaches are based on finite dimensional approximations with the assumption of  an explicitly given moving boundary dynamics \cite{Daraoui2010,Armaou01}. 
 For control objectives, infinite-dimensional approaches have been used for stabilization of  the temperature profile and the moving interface of a 1D Stefan problem, such as enthalpy-based feedback~\cite{petrus2012} and geometric control~\cite{maidi2014}.  These works designed control laws ensuring the asymptotical stability of the closed-loop system in the ${L}_2$ norm. However, the results in \cite{maidi2014} are established based on the assumptions on the liquid temperature being greater than the melting temperature, which must be ensured by showing the positivity of the boundary heat input. 
 
Recently, boundary feedback controllers for the Stefan problem have been designed via a ``backstepping transformation" \cite{krstic2008boundary,krstic2009delay,andrew2004} which has been used for many other classes of infinite-dimensional systems. For instance, \cite{Shumon16} designed a state feedback control law by introducing a nonlinear backstepping transformation for moving boundary PDE, which achieved the exponentially stabilization of the closed-loop system in the ${\mathcal H}_1$ norm without imposing any {\em a priori} assumption. Based on the technique, \cite{Shumon16CDC} designed an observer-based output feedback control law for the Stefan problem, \cite{Shumon19journal} extended the results in \cite{Shumon16, Shumon16CDC} by studying the robustness with respect to the physical parameters and developed an analogous design with Dirichlet boundary actuation, \cite{Shumon2017ACC} designed a state feedback control for the Stefan problem under the material's convection, and \cite{koga2017CDC,koga_2019delay} developed a delay-compensated control for the Stefan problem with proving the robustness to the mismatch between the actuator's delay and the compensated delay.  

While all the aforementioned results are based on the one-phase Stefan problem which neglects the cooling heat caused by the solid phase, an analysis on the system incorporating the cooling heat at the liquid-solid interface has not been established. The one-phase Stefan problem with a prescribed heat flux at the interface was studied in \cite{Sherman67}. The author proved the existence and uniqueness of the solution with a prescribed heat input at the fixed boundary by verifying positivity conditions on the interface position and temperature profile using a similar technique as in \cite{Friedman59}. 

Regarding the added heat flux at the interface as the heat loss induced by the remaining other phase dynamics, it is reasonable to treat the prescribed heat flux as the disturbance of the system. The norm estimate of systems with a disturbance is often analyzed in terms of ISS \cite{Sontag08}, which serves as a criteria for the robustness of the controller or observer design \cite{arcak01}. Recently, the ISS for infinite dimensional systems with respect to the boundary disturbance was developed in \cite{Karafyllis16,Karafyllis17,Karafyllis-issbook} using the spectral decomposition of the solution of linear parabolic PDEs in one dimensional spatial coordinate with Strum-Liouville operators. An analogous result for the diffusion equations with a radial coordinate in $n$-dimensional balls is shown in \cite{Leobardo2018} with proposing an application to robust observer design for battery management systems \cite{Moura2014}.  

The present paper revisits the result of our conference paper \cite{koga2018ISS} which proved ISS of the one-phase Stefan problem under the control law proposed in \cite{Shumon16} with respect to the heat loss at the interface, and extends the ISS result to the two-phase Stefan problem with an unknown heat loss at the boundary of the solid phase under the control law developed in \cite{koga2018CDC}. First, we consider a prescribed open-loop control of the one-phase Stefan problem in which the solution of the system is equivalent with the proposed closed-loop control. Using the result of \cite{Sherman67}, the well-posedness of the solution and the positivity conditions are proved. Next, we apply the closed-loop control through the backstepping transformation as in \cite{Shumon16}. The well-posedness of the closed-loop solution and the positivity conditions for the model to be valid are ensured by showing the differential equation of the control law. The associated target system has the disturbance at the interface dynamics due to the heat loss. An estimate of the ${ L}_{2}$ norm of the closed-loop system is developed in the sense of ISS through the analysis of the target system by Lyapunov method. The similar procedure is performed to extend the ISS result to the two-phase Stefan problem with the unknown heat loss at the boundary of the solid phase. 

This paper is organized as follows. The one-phase Stefan problem with heat loss at the interface is presented  in Section \ref{model}. The well-posedness of the system with an open-loop controller is proved in Section \ref{open}. Section \ref{sec:closed} introduces our first main result for the ISS of the closed-loop system under the designed feedback control law with providing the proof. The extension of the result to the two-phase Stefan problem is addressed in Section \ref{sec:2ph}. Supportive numerical simulations are provided in Section \ref{simulation}. The paper ends with some final remarks in Section \ref{conclusion}.

\section{Description of the Physical Process}\label{model}
\begin{figure}[t]
\centering
\includegraphics[width=2.5in]{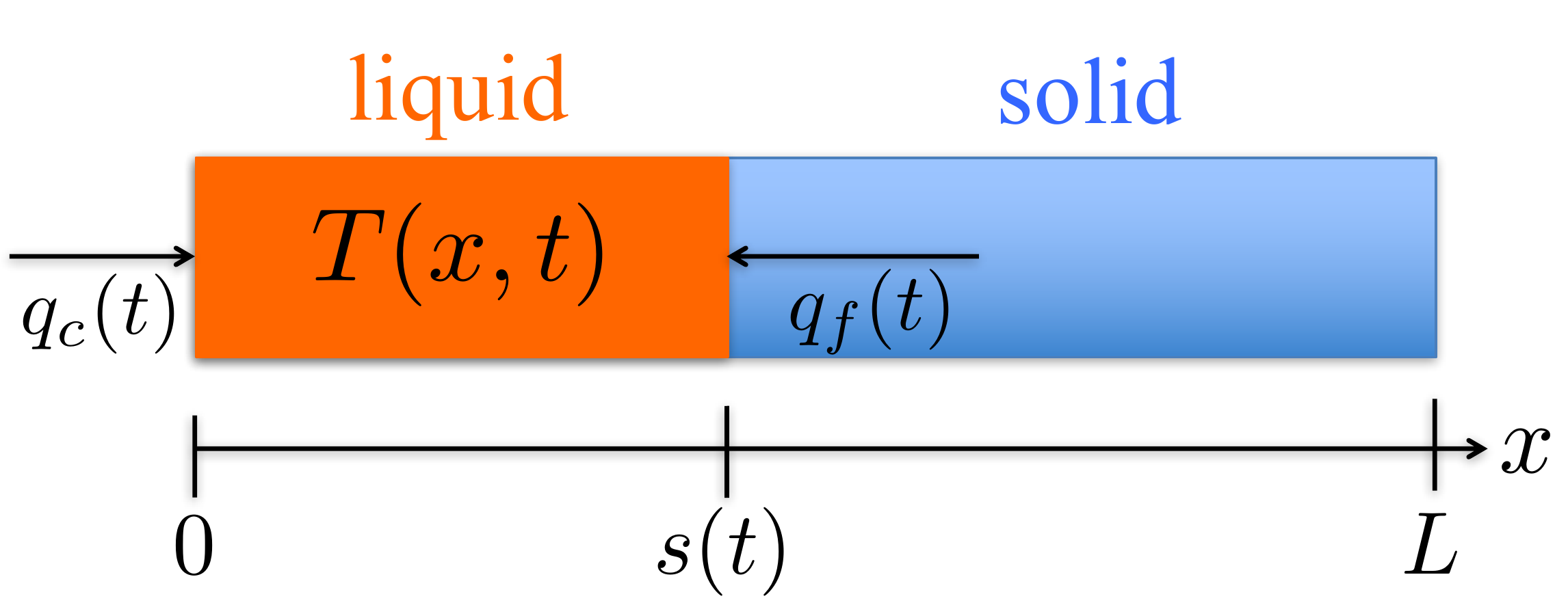}\\
\caption{Schematic of the one-phase Stefan problem.}
\label{fig:stefan}
\end{figure}

Consider a physical model which describes the melting or solidification mechanism in a pure one-component material of length $L$ in one dimension. In order to mathematically describe the position at which phase transition occurs, we divide the domain $[0, L]$ into two time-varying sub-domains, namely,  the interval $[0,s(t)]$ which contains the liquid phase, and the interval $[s(t),L]$ that contains the solid  phase. A heat flux enters the material through the boundary at $x=0$ (the fixed boundary of the liquid phase) which affects the liquid-solid interface dynamics through heat propagation in liquid phase. In addition, due to the cooling effect from the solid phase, there is a heat loss at the interface position $x  = s(t)$. As a consequence, the heat equation alone does not provide a complete description of the phase transition and must be coupled with the dynamics that describes the moving boundary. This configuration is shown in Fig.~\ref{fig:stefan}.

Assuming that the temperature in the liquid phase is not lower than the melting temperature  of the material $T_{{\mathrm m}}$, the energy conservation and heat conduction laws yield the heat equation of  the liquid phase as follows
\begin{align}\label{eq:stefanPDE}
T_t(x,t)&=\alpha T_{xx}(x,t), \quad  \alpha :=\fr{k}{\rho C_p}, \quad  0\leq x\leq s(t), 
\end{align}
with the boundary conditions
\begin{align}\label{eq:stefancontrol}
-k T_x(0,t)&=q_{\rm c}(t), \\ \label{eq:stefanBC}
T(s(t),t)&=T_{{\rm m}},
\end{align}
and the initial values
\begin{align}\label{eq:stefanIC}
T(x,0)=T_0(x), \quad s(0) = s_0
\end{align}
where $T(x,t)$,  ${q}_c(t)$,  $\rho$, $C_p$, and $k$ are the distributed temperature of the liquid phase, the manipulated heat flux, the liquid density, the liquid heat capacity, and the liquid heat conductivity, respectively. The first problem we consider in this paper assumes that the heat loss at the interface caused by the solid phase temperature dynamics is assumed to be a prescribed function in time, denoted as $q_{\rm f}(t)$. The local energy balance at the liquid-solid interface $x=s(t)$ yields
\begin{align}\label{eq:stefanODE}
\rho \Delta H^* \dot{s}(t)=- k T_x(s(t),t) - q_{\rm f}(t), 
\end{align}
where $\Delta H^*$ represents the latent heat of fusion. In \eqref{eq:stefanODE}, the left hand side represents the latent heat, and the first and second term of the right hand side represent the heat flux by the liquid phase and the heat loss caused by the solid phase, respectively. As the governing equations \eqref{eq:stefanPDE}--\eqref{eq:stefanODE} suffice to determine the dynamics of the states $(T,s)$, the temperature in the solid phase is not needed to be modeled. 

\begin{remark}As the moving interface  $s(t)$ depends on the temperature, the problem defined in  \eqref{eq:stefanPDE}--\eqref{eq:stefanODE}  is nonlinear.\end{remark}

There are two underlying assumptions to validate the model \eqref{eq:stefanPDE}-\eqref{eq:stefanODE}. First, the liquid phase is not frozen to solid from the boundary $x=0$. This condition is ensured if the liquid temperature $T(x,t)$ is greater than the melting temperature. Second, the material is not completely melt or frozen to single phase through the disappearance of the other phase. This condition is guaranteed if the interface position remains inside the material's domain. In addition, these conditions are also required for the well-posedness (existence and uniqueness) of the solution in this model. Taking into account of these model validity conditions, we emphasize the following remark. 
\begin{remark}
To maintain the model \eqref{eq:stefanPDE}-\eqref{eq:stefanODE} to be physically validated, the following conditions must hold: 
\begin{align}\label{temp-valid}
T(x,t) \geq& T_{{\rm m}}, \quad  \forall x\in(0,s(t)), \quad \forall t>0, \\
\label{int-valid}s(t)>&0, \quad \forall t>0. 
\end{align}
\end{remark}

Based on the above conditions, we impose the following assumption on the initial data. 
\begin{assum}\label{initial} 
$s_0>0$, $T_0(x) \geq T_{{\rm m}}$ for all $x \in [0,s_0]$, and $T_0(x) $ is continuously differentiable in $x\in[0,s_0]$.
 \end{assum}
\begin{lem}\label{lem1}
Under Assumption \ref{initial}, if the boundary input keeps generating positive heat, i.e., $q_{{\rm c}}(t) \geq 0, $ for all $t>0$, then the condition \eqref{temp-valid} holds. 
\end{lem}

The proof of Lemma \ref{lem1} is established by maximum principle and Hopf's lemma as shown in \cite{Gupta03}. For the heat flux at the interface, the following assumptions are imposed. 
\begin{assum}\label{disturbed}
The heat loss remains non-negative, bounded, and continuous for all $ t\geq 0$, and the total energy is also bounded, i.e., 
\begin{align}
q_{\rm f}(t)\geq& 0, \quad \forall t>0, \\
\label{bound}\exists M>&0, \quad \textrm{s.t.} \quad \int_0^{\infty} q_{\rm f}(t) {\rm d}t < M. 
\end{align}
\end{assum}

The steady-state solution $(T_{{\rm eq}}(x), s_{{\rm eq}})$ of the system \eqref{eq:stefanPDE}-\eqref{eq:stefanODE} with zero manipulating heat flux $q_{\rm c}(t)=0$ and zero heat loss at the interface $q_{\rm f}(t) = 0$ yields a uniform melting temperature $T_{{\rm eq}}(x) = T_{{\mathrm m}}$ and a constant interface position given by the initial data. In \cite{Shumon16}, the authors developed the exponential stabilization of the interface position $s(t)$ at a desired reference setpoint $s_{{\mathrm r}}$ with zero heat loss $q_{{\rm f}}(t)=0$ through the following state feedback control design of $q_{{\rm c}}(t)$: 
\begin{align}\label{control}
q_{\rm c}(t) = - c \left( \fr{k}{\alpha } \int_0^{s(t)} (T(x,t)-T_{{\rm m}}) {\rm d}x + \fr{k}{\beta} (s(t)-s_{{\rm r}}) \right),   
\end{align}
where $c>0$ is the control gain which can be chosen by user. To maintain the positivity of the heat input $q_{{\rm c}}(t)$ as stated in Lemma \ref{lem1}, the following restriction on the setpoint is imposed. 
\begin{assum}\label{ass:setpoint}
The setpoint is chosen to verify
\begin{align}\label{eq:setpoint}
s_{{\rm r}} > s_0 + \fr{\beta}{\alpha} \int_0^{s_0} (T_0(x)-T_{{\rm m}}) {\rm d}x , 
\end{align}
where $\beta := \fr{k}{\rho \Delta H^*}$. 
\end{assum}

Finally, we impose the following condition of the control gain. 
\begin{assum}\label{ISS-cond}
The control gain $c$ is chosen sufficiently large to satisfy $  c   > \fr{\beta}{k s_{{\rm r}}} \bar{q}_{\rm f}$, where $\bar{q}_{\rm f} := \sup_{0\leq t \leq \infty} \left\{q_{\rm f}(t)\right\}$. 
\end{assum}

 In this paper, we prove the ISS of the reference error system with the control design \eqref{control} with respect to the heat loss at the interface by studying the norm estimate through Sections \ref{open}--\ref{sec:closed}. \\

\section{Open-Loop System and Analysis} \label{open}
The well-posedness of the one-phase Stefan problem with heat flux at the interface was developed in \cite{Sherman67} with a prescribed open-loop heat input for $q_{{\rm c}}(t)$. To apply the result, in this section we focus on an open-loop control which has an equivalent solution as the closed-loop control introduced later, and we prove the well-posedness of the open-loop system. 
\begin{lem}\label{posed}
Let Assumptions \ref{initial}-\ref{ass:setpoint} hold. For any $\bar t \leq \sigma $ where $0<\sigma \leq \infty$, there is a unique solution of the system \eqref{eq:stefanPDE}--\eqref{eq:stefanODE} with the open-loop control 
\begin{align}\label{cont-exp}
q_{\rm c}(t) = q_{0} e^{-ct} + c \int_0^t e^{-c(t-\tau)} q_{\rm f}(\tau) {\rm d}\tau , 
\end{align}
for $0<t<\bar t$, where
\begin{align}\label{q0}
q_{0} = - c \left( \fr{k}{\alpha } \int_0^{s_0} (T_0(x)-T_{{\rm m}}) {\rm d}x + \fr{k}{\beta} (s_0-s_{{\rm r}}) \right). 
\end{align}
If $\sigma \neq \infty$, then $s(\sigma) = 0$. 
\end{lem}
\begin{pf}
Since Assumption \ref{setpoint} leads to $q_{0}>0$, the open-loop controller \eqref{cont-exp} remains positive and continuous for all $ t>0$ by Assumption \ref{disturbed}. Hence, applying the theorem in \cite{Sherman67} (page 3) proves Lemma \ref{posed}. 
\end{pf}

\begin{lem}\label{global}
Under Assumption \ref{ISS-cond}, Lemma \ref{posed} holds globally, i.e., $\sigma=\infty$. 
\end{lem}
\begin{pf}
We prove by contradiction. Suppose there exists $0<\sigma^*<\infty$ such that $s(\sigma^*)=0$. Let $E(t)$ be an internal energy of the system defined by 
\begin{align}\label{Et}
E(t) = \fr{k}{\alpha} \int_0^{s(t)} \left(T(x,t) - T_{{\rm m}} \right) {\rm d}x + \fr{k}{\beta} s(t) . 
\end{align}
with $E(\sigma^*) = 0$ by the imposed assumption. Taking the time derivative of \eqref{Et} yields the energy conservation 
\begin{align} \label{energy-1ph} 
\dot{E}(t) = q_{\rm c}(t) - q_{\rm f}(t) . 
\end{align} 
In addition, the time derivative of \eqref{cont-exp} yields $\dot{q}_{\rm c}(t) = -c \left(q_{\rm c}(t) - q_{\rm f}(t) \right)$. Combining these two and taking integration on both sides gives
\begin{align}\label{EtE0}
E(t) = E(0)- \fr{1}{c} (q_{\rm c}(t)-q_{\rm c}(0)) . 
\end{align}
By \eqref{q0} and \eqref{Et}, we get $q_{0} = - c (E(0) - \fr{k}{\beta} s_{{\rm r}})$. Substituting this and \eqref{cont-exp} into \eqref{EtE0}, we have 
\begin{align}
E(t) = e^{-ct} \left[ E(0) + \fr{k s_{{\rm r}}}{\beta} (e^{ct} - 1) -  \int_0^{t} e^{c \tau } q_{\rm f}(\tau ) {\rm d}\tau \right] . 
\end{align}
Let $f(t)$ be a function in time defined by
\begin{align}\label{time-der-f}
f(t) = E(0) + \fr{k s_{{\rm r}}}{\beta} (e^{ct} - 1) -  \int_0^{t} e^{c \tau } q_{\rm f}(\tau ) {\rm d}\tau. 
\end{align}
Since $E(t) = e^{-ct} f(t)$, we can see that $E(t)>0$ for all $t>0$ if and only if $f(t)>0$ for all $t>0$. By \eqref{time-der-f}, we have $f(0) =  E(0)>0$. Taking the time derivative of \eqref{time-der-f} yields 
\begin{align}\label{fprime}
f'(t) = e^{ct} \left( \fr{k  s_{{\rm r}}c }{\beta} - q_{\rm f} (t) \right) . 
\end{align}
By Assumption \ref{ISS-cond}, \eqref{fprime} leads to $f'( t)>0$ for all $ t>0$. Therefore, $f(t)>0$ for all $t>0$, and we conclude $E(t)>0$ for all $t>0$ which contradicts with the imposed assumption $s(\sigma^*)=0$ where $\sigma^* \neq \infty$. Hence, Lemma \ref{posed} holds globally.  
\end{pf}
\begin{coro} \label{properties}
The open-loop solution satisfies the model validity conditions \eqref{temp-valid} \eqref{int-valid} for all $t>0$. 
\end{coro}

As proven in \cite{Sherman67}, the global well-posedness shown in Lemma \ref{global} verifies the conditions  \eqref{temp-valid} \eqref{int-valid} for all $ t>0$ which derives Corollary \ref{properties}. 

\section{ISS for One-Phase Stefan Problem}\label{sec:closed}
While the open-loop input \eqref{cont-exp} ensures the well-posed solution of the system, the analysis does not enable to prove the ISS property of the norm estimate. In addition, the open-loop design requires the heat loss at the interface $q_{\rm f}(t)$. However, measuring $q_{\rm f}(t)$ is not practically doable. In this section, we show that the closed-loop solution with the control law proposed in \cite{Shumon16} is equivalent to the open-loop solution introduced in Section \ref{open}. The controller is feedback design of liquid temperature profile and the interface position $(T(x,t), s(t))$. The heat loss $q_{{\rm f}}(t)$ is regarded as a disturbance, and the norm estimate of the reference error is derived in a sense of input-to-state stability. 

Our first main result is stated in the following theorem. 

\begin{thm}\label{theo-1}
Under Assumptions \ref{initial}-\ref{ISS-cond}, the closed-loop system consisting of \eqref{eq:stefanPDE}--\eqref{eq:stefanODE} with the control law \eqref{control} 
satisfies the model validity conditions \eqref{temp-valid} and \eqref{int-valid}, and is ISS with respect to the heat loss $q_{{\rm f}}(t)$ at the interface, i.e., there exist a class-$\mathcal{KL}$ function $\zeta$ and a class-${\mathcal K}$ function $\eta$ such that the following estimate holds:
\begin{align}\label{h1}
\Psi(t) \leq \zeta(\Psi(0),t) + \eta \left( \sup_{\tau \in [0,t]} |q_{{\rm f}}(\tau)| \right),
\end{align}
for all $t\geq0$, in the $L_2$ norm \\
$\Psi(t) = \left( \int_{0}^{s(t)} \left(T(x,t)-T_{{\rm m}}\right)^2 {\rm d}x +(s(t)-s_{r})^2 \right)^{\frac{1}{2}}$. Moreover, there exist positive constants $M_1>0$ and $M_2>0$ such that the explicit functions of $\zeta$ and $\eta$ are given by $\zeta(\Psi(0),t) =M_1 \Psi(0) e^{- \lambda t}$, $\eta( \sup_{\tau \in [0,t]} |q_{{\rm f}}(\tau)|)= M_2  \sup_{\tau \in [0,t]} |q_{{\rm f}}(\tau)| $, where $\lambda =\fr{1}{32} \min \left\{ \fr{\alpha}{s_{{\rm r}}^2}, c \right\} $, which ensures the exponentially ISS.  
\end{thm}

The proof of Theorem \ref{theo-1} is established through the remainder of this section. 

\subsection{Reference error system}\label{Ref}
Let $u(x,t)$ and $X(t)$ be reference error variables defined by 
\begin{align} 
u(x,t):= &T(x,t)-T_{{\rm m}}, \\
X(t):= &s(t)-s_{r}. 
\end{align} 
Then, the system \eqref{eq:stefanPDE}--\eqref{eq:stefanODE} is rewritten as
\begin{align}\label{u-sys1}
u_{t}(x,t) =&\alpha u_{xx}(x,t),\\
\label{u-sys2}u_x(0,t) =& - \fr{q_{{\rm c}}(t)}{k},\\
\label{u-sys3}u(s(t),t) =&0,\\
\label{u-sys4}\dot{X}(t) =&-\beta u_x(s(t),t) - d(t),
\end{align}
where $d(t) = \fr{\beta}{k} q_{\rm f}(t)$. The controller is designed to stabilize $(u,X)$-system at the origin for $d(t)=0$. 
\subsection{Backstepping transformation} 
Introduce the following backstepping transformation 
\begin{align}\label{eq:DBST}
w(x,t)=&u(x,t)-\frac{\beta}{\alpha} \int_{x}^{s(t)} \phi (x-y)u(y,t) {\rm d}y \notag\\
&-\phi(x-s(t)) X(t),
\end{align}
which maps into 
\begin{align}\label{tarPDE}
w_t(x,t)=&\alpha w_{xx}(x,t)+ \dot{s}(t) \phi'(x-s(t))X(t) \notag\\
&+\phi(x-s(t)) d(t), \\
\label{tarBC2} w_x(0,t) =&  \fr{\beta}{\alpha} \phi(0) u(0), \\
\label{tarBC1} w(s(t),t) =& \ep X(t), \\
\label{tarODE}\dot{X}(t)=&-cX(t)-\beta w_x(s(t),t) - d(t).
\end{align}
The objective of the transformation \eqref{eq:DBST} is to add a stabilizing term $-c X(t)$ in \eqref{tarODE} of the target $(w,X)$-system which is easier to prove the stability for $d(t)=0$ than $(u,X)$-system. 
By taking the derivative of \eqref {eq:DBST} with respect to $t$ and $x$ respectively along the solution of \eqref{u-sys1}-\eqref{u-sys4}, to satisfy \eqref{tarPDE}, \eqref{tarBC1}, \eqref{tarODE}, we derive the conditions on the gain kernel solution which yields the solution as 
\begin{align}
\label{kernel}\phi(x) =& \frac{c}{\beta} x- \ep. 
\end{align}
By matching the transformation \eqref{eq:DBST} with the boundary condition \eqref{tarBC2}, the control law is derived as 
\begin{align} \label{qcfb}
q_{{\rm c}}(t) = -c \left(\frac{k}{\alpha} \int_{0}^{s(t)} u(y,t) {\rm d}y + \fr{k}{\beta} X(t) \right) . 
\end{align}
Rewriting \eqref{qcfb} by $T(x,t)$ and $s(t)$ yields \eqref{control}. 
\subsection{Inverse transformation} 
Consider the following inverse transformation 
\begin{align}\label{inv-trans}
u(x,t)=&w(x,t)-\frac{\beta}{\alpha} \int_{x}^{s(t)} \psi (x-y)w(y,t) {\rm d}y \notag\\
&-\psi(x-s(t)) X(t). 
\end{align}
Taking the derivatives of \eqref{inv-trans} in $x$ and $t$ along \eqref{tarPDE}-\eqref{tarODE}, to match with \eqref{u-sys1}-\eqref{u-sys4}, we obtain the gain kernel solution as
\begin{align}\label{inv-gain}
\psi(x) =& e^{ r x } \left( p_1 \sin\left( \omega x \right) + \ep \cos\left( \omega x \right) \right) , 
\end{align}
where $r = \fr{\beta \varepsilon}{2 \alpha}$, $ \omega = \sqrt{\fr{4 \alpha c - (\varepsilon\beta)^2 }{4 \alpha^2 } }$, $p_1 = - \fr{1}{2 \alpha \beta \omega} \left( 2 \alpha c - (\varepsilon \beta )^2 \right) $, 
and $0<\ep<2 \fr{\sqrt{\alpha c}}{\beta}$ is to be chosen later. Finally, using the inverse transformation, the boundary condition \eqref{tarBC2} is rewritten as 
\begin{align}
\label{eq:tarBC2} w_x(0,t) =&  - \frac{\beta}{\alpha}\varepsilon \left[  w(0,t)-\frac{\beta}{\alpha} \int_{0}^{s(t)} \psi (-y)w(y,t) {\rm d}y \right. \notag\\
&\left. -\psi(-s(t)) X(t) \right]. 
\end{align}
In other words, the target $(w,X)$-system is described by \eqref{tarPDE}, \eqref{tarBC1}, \eqref{tarODE}, and \eqref{eq:tarBC2}. Note that the boundary condition \eqref{tarBC1} and the kernel function \eqref{kernel} are modified from the one in \cite{Shumon16}, while the control design \eqref{qcfb} is equivalent.  The target system derived in \cite{Shumon16} requires ${\mathcal H}_1$-norm analysis for stabilty proof. However, with the prescribed heat loss at the interface, ${\mathcal H}_1$-norm analysis fails to show the stability due to the non-monotonic moving boundary dynamics. The modification of the boundary condition \eqref{tarBC1} enables to prove the stability in ${ L}_2$ norm as shown later. 
\subsection{Analysis of closed-loop system} \label{sec:wellposed}
Here, we prove the well-posedness of the closed-loop solution and the positivity conditions of the state variables. 
\begin{lem}\label{closed-posed}
Under Assumptions \ref{initial}-\ref{ISS-cond}, the closed-loop system of \eqref{eq:stefanPDE}--\eqref{eq:stefanODE} with the control law \eqref{control} has a unique solution which is equivalent to the open-loop solution of \eqref{eq:stefanPDE}--\eqref{eq:stefanODE} with the control law \eqref{cont-exp} for all $t>0$. 
\end{lem}
\begin{pf} 
Taking the time derivative of the control law \eqref{control} along with the energy conservation leads to the following differential equation 
\begin{align} 
\dot{q}_{\rm c} (t) = - c q_{\rm c}(t) + c q_{\rm f}(t), 
\end{align} 
which has the same explicit solution as the open-loop control \eqref{cont-exp}. Hence, the closed-loop solution is equivalent to to the open-loop solution with \eqref{cont-exp}. Since the open-loop system has a unique solution as shown in Lemma \ref{global}, the closed-loop solution has a unique solution as well.   
\end{pf}

\begin{lem}\label{closed-validity}
The closed-loop solution verifies the following properties: 
\begin{align}\label{qcpositive}
q_{\rm c}(t)>&0, \\
\label{upositive} u(x,t)>&0, \quad u_{x}(s(t),t)<0, \\
 s(t)>&0. 
\end{align}
\end{lem}
\begin{pf}
Applying Corollary \ref{properties} and Lemma \ref{closed-posed} gives the proof of Lemma \ref{closed-validity}. 
\end{pf}
\begin{lem}
The closed-loop system has the following property:
\begin{align}\label{overshoot}
0<s(t)< s_{{\rm r}} . 
\end{align}
\end{lem}
\begin{pf}
Applying \eqref{qcpositive} and \eqref{upositive} to \eqref{qcfb}, the condition \eqref{overshoot} is verified.
\end{pf}

  \subsection{Stability analysis} \label{sec:stability}
To conclude the ISS of the original system, first we show the ISS of the target system \eqref{tarPDE}, \eqref{tarBC1}, \eqref{tarODE}, and \eqref{eq:tarBC2} with respect to the disturbance $d(t)$. We consider the following functional
\begin{align}\label{lyap}
V(t) = \fr{1}{2\alpha } || w||^2 + \fr{\varepsilon}{2\beta } X(t)^2, 
\end{align}
where $||w||$ denotes $L_2$ norm defined by $||w|| = \left( \int_0^{s(t)} w(x,t)^2 {\rm d}x \right)^{1/2} $. Taking the time derivative of \eqref{lyap} along with the solution of   \eqref{tarPDE}--\eqref{tarODE}, \eqref{eq:tarBC2}, we have 
\begin{align}\label{Vdot}
\dot{V} (t)=&  - || w_{x}||^2  - \fr{\varepsilon}{\beta}cX(t)^2 +\frac{\beta}{\alpha}\varepsilon w(0,t)^2 \notag\\
& -\frac{\beta}{\alpha}\varepsilon w(0,t)\left[ \frac{\beta}{\alpha} \int_{0}^{s(t)} \psi (-y)w(y,t) {\rm d}y \right. \notag\\
&\left. +\psi(-s(t)) X(t) \right] \notag\\
&- \fr{\varepsilon}{\beta } X(t) d(t)+ \fr{1}{\alpha} \int_0^{s(t)} \phi(x-s(t))w(x,t) {\rm d}x d(t)  \notag\\
&+ \fr{\dot{s}(t)}{\alpha} \left( \fr{\ep^2}{2 } X(t)^2+  \fr{ c}{\beta} \int_0^{s(t)} w(x,t) {\rm d}x X(t)  \right) . 
\end{align}
Applying Young's inequality to the second and third lines of \eqref{Vdot} twice, we get
 \begin{align}\label{dyoung}
 & - w(0,t)\left[\frac{\beta}{\alpha} \int_{0}^{s(t)} \psi (-y)w(y,t) dy+\psi(-s(t)) X(t) \right] \notag\\
  \leq& \fr{1}{2}w(0,t)^2 + \frac{\beta^2}{\alpha^2 \gamma_1} \left(\int_{0}^{s(t)} \psi (-y)w(y,t) {\rm d}y \right)^2 \notag\\
  &+ \gamma_1 \left(\psi(-s(t)) X(t)\right)^2. 
  \end{align}
  Again applying Young's inequality to the forth line of \eqref{Vdot}, 
  \begin{align}\label{dyoung2}
&- \fr{\varepsilon}{\beta } X(t) d(t)+ \fr{1}{\alpha} \int_0^{s(t)} \phi(x-s(t))w(x,t) {\rm d}x d(t) \notag\\
\leq& \fr{1}{2\gm_2} \left( \fr{\varepsilon}{\beta } X(t) \right)^2 + \fr{\left( \gm_2 + \gm_3 \right)}{2} d(t)^2 \notag\\
& + \fr{1}{2\alpha^2 \gm_3}\left( \int_0^{s(t)} \phi(x-s(t))w(x,t) {\rm d}x \right)^2, 
\end{align}
where $\gm_{i}>0$ for $i=\{1,2,3\}$. Applying \eqref{dyoung}-\eqref{dyoung2} and Cauchy Schwarz, Poincare, and Agmon's inequalities to \eqref{Vdot} with choosing $\gm_1 = \fr{1}{8}$, $\gm_2=\fr{2 \ep}{\beta c}$, and $\gm_3 = \frac{4 s_{{\rm r}}^3}{\alpha^2} \left( \fr{c s_{{\rm r}}}{\beta} + \ep \right)^2$,  we have 
\begin{align}\label{Vdot3}
\dot{V} (t) \leq &  - \left(\fr{1}{2} - \fr{2 \beta s_{{\rm r}}}{\alpha} \left( \frac{ 64 c s_{{\rm r}}^2}{\alpha } + 3 \right)\varepsilon\right) || w_{x}||^2 \notag\\
& - \varepsilon \left(\fr{c}{8\beta} + g(\ep) \right)X(t)^2 + \Gamma d(t)^2    \notag\\
&+ \fr{|\dot{s}(t)|}{2\alpha} \left( \ep^2 X(t)^2 +  \fr{ 2c}{\beta} \left|\int_0^{s(t)} w(x,t) {\rm d}x X(t)  \right|\right) , 
\end{align}
where $\Gamma =  \fr{\left( \gm_2 + \gm_3 \right)}{2} $, and $g(\ep) = \fr{c}{8\beta} - \fr{\ep}{4 s_{{\rm r}}} -  \fr{ \beta }{\alpha} \left( \frac{ 64 c s_{{\rm r}}^2}{\alpha } + 3 \right) \ep^2 $. Since $g(0) = \fr{c}{8\beta }>0$ and $g'(\ep) = - \fr{1}{4 s_{{\rm r}}} - \frac{2\beta \ep }{\alpha}  \left( \frac{ 64 c s_{{\rm r}}^2}{\alpha } + 3 \right) <0$ for all $\ep>0$, there exists $\ep^*$ such that $g(\ep)>0$ for all $\ep \in (0, \ep^*)$ and $g(\ep^*)=0$. Thus, setting $\ep < \min\left\{ \ep^*, \fr{\alpha}{8 \beta s_{{\rm r}} \left( \frac{ 64 c s_{{\rm r}}^2}{\alpha } + 3 \right)} \right\}$ , 
the inequality \eqref{Vdot3} leads to
\begin{align}\label{dotV}
\dot{V} (t) \leq &  - b V (t) + \Gamma d(t)^2    \notag\\
&\hspace{-0mm}+ \fr{|\dot{s}(t)|}{2\alpha} \left( \ep^2 X(t)^2 +  \fr{2 c}{\beta} \left|\int_0^{s(t)} w(x,t) {\rm d}x X(t)  \right|\right) , 
\end{align}
where $b =\fr{1}{8} \min \left\{ \fr{\alpha}{s_{{\rm r}}^2}, c \right\} $. 
Since $u_{x}(s(t),t)<0$ and $d(t)>0$, we have 
\begin{align}\label{sdot}
|\dot{s}(t) | \leq - \beta u_{x}(s(t),t) +  d(t). 
\end{align}
 Let $z(t)$ be defined by 
 \begin{align}\label{zdef}
z(t) := s(t) + 2 \int_0^{t} d(\tau) {\rm d}\tau . 
\end{align}
The time derivative of \eqref{zdef} becomes $\dot{z}(t) = - \beta u_{x}(s(t),t) + d(t)$. In addition, applying \eqref{bound} and \eqref{overshoot} to \eqref{zdef} leads to 
\begin{align} \label{zbound} 
0< z(0) < z(t) < s_{{\rm r}} + \fr{2 \beta M}{k} := \bar z . 
\end{align} 
Applying \eqref{sdot} and Young's inequality to \eqref{dotV}, we have 
\begin{align}\label{dotV2}
\dot{V}(t) \leq &  - b V (t) + \Gamma d(t)^2   + a\dot{z}(t) V(t) , 
\end{align}
where $a = \fr{2\beta \ep }{\alpha} \max \left\{1,\fr{\alpha c^2 s_{{\rm r}}}{2\beta^3 \ep^3} \right\}$. Consider the following functional 
\begin{align}\label{Wdef}
W (t)= V(t) e^{ - a z(t)} . 
\end{align}
Taking the time derivative of \eqref{Wdef} with the help of \eqref{dotV2} and applying \eqref{zbound}, we deduce 
\begin{align}\label{Wdot}
\dot{W} (t) \leq - b W (t) + \Gamma d(t)^2 . 
\end{align}
Since \eqref{Wdot} leads to the statement that either $\dot{W} (t)\leq - \frac{b}{2} W(t)$ or $W(t) \leq \frac{2}{b} \Gamma d(t)^2$ is true, following the procedure in \cite{Sontag08} (proof of Theorem 5 in Section 3.3), one can derive 
\begin{align} \label{Wineq} 
W(t) \leq W(0) e^{- \frac{b}{2} t} + \frac{2}{b} \Gamma \sup_{\tau \in [0.t]} d(\tau)^2. 
\end{align} 
By \eqref{Wdef} and \eqref{Wineq} and applying \eqref{zbound}, finally we obtain the following estimate on the $L_2$ norm of the target system
\begin{align}
V(t) \leq & V(0) e^{a \bar{z}} e^{- \frac{b}{2}t} + \frac{2}{b} \Gamma  e^{a\bar z}  \sup_{\tau \in [0.t]} d(\tau)^2 . 
\end{align}
By the invertibility of the transformation, the norm equivalence of the target $(w,X)$-system to the $(u,X)$-system holds, and hence the norm estimate on the original $(T,s)$-system is derived, which completes the proof of Theorem \ref{theo-1}. 

\section{ISS for Two-Phase Stefan Problem} \label{sec:2ph} 
\begin{figure}[t]
\centering
\includegraphics[width=2.5in]{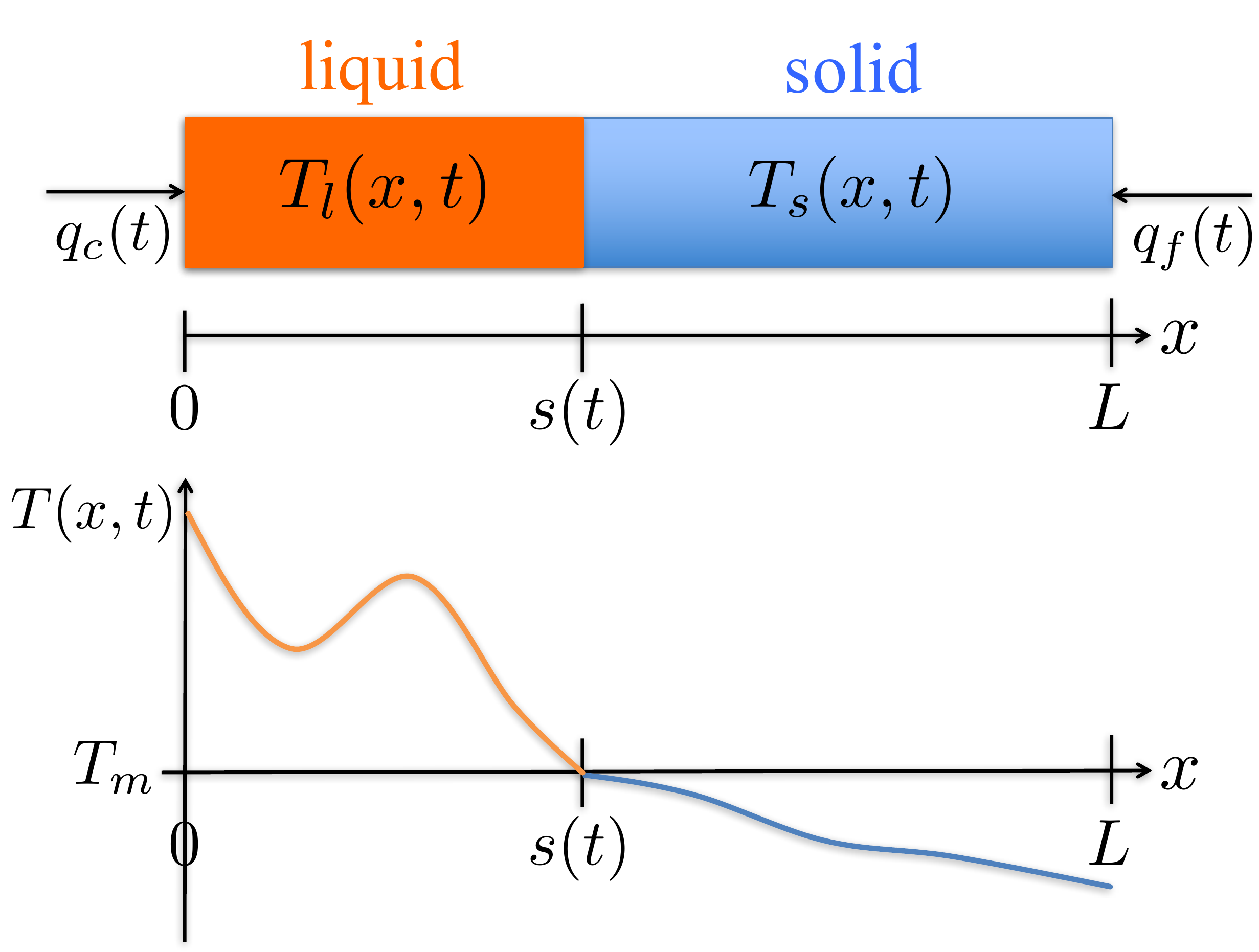}\\
\caption{Schematic of the two-phase Stefan problem.}
\label{fig:2phstefan}
\end{figure}
In this section, we extend the results we have established up to the last section to the "two-phase" Stefan problem, where the heat loss at the interface is precisely modeled by the temperature dynamics in the solid phase, following the work in \cite{koga2018CDC}. An unknown heat loss is then accounted not at the interface, but at the boundary of the solid phase. This configuration is depicted in Fig. \ref{fig:2phstefan}. 
\subsection{Problem statement} 
The governing equations are descried by the following coupled PDE-ODE-PDE system: 
\begin{align}\label{sys1}
 \fr{\pa T_{{\rm l}}}{\pa t}(x,t) =&\alpha_{{\rm l}}  \fr{\pa^2 T_{{\rm l}}}{\pa x^2}(x,t), \quad 0<x<s(t),\\
 \fr{\pa T_{{\rm l}}}{\pa x}(0,t) =& -q_{{\rm c}}(t)/k_{{\rm l}}, \quad T_{{\rm l}}(s(t),t) =T_{{\rm m}},\\
\label{sys3} \fr{\pa T_{{\rm s}}}{\pa t}(x,t) =&\alpha_{{\rm s}}  \fr{\pa^2 T_{{\rm s}}}{\pa x^2}(x,t), \quad s(t)<x<L, \\
\label{sys4} \fr{\pa T_{{\rm s}}}{\pa x}(L,t) =& - q_{{\rm f}}(t)/k_{{\rm s}},\quad T_{{\rm s}}(s(t),t) =T_{{\rm m}}, \\
\label{sys5} \gm \dot{s}(t) =& - k_{{\rm l}} \fr{\pa T_{{\rm l}}}{\pa x}(s(t),t)+k_{{\rm s}} \fr{\pa T_{{\rm s}}}{\pa x}(s(t),t),
\end{align}
where $\gamma = \rho_{{\rm l}} \Delta H^*$, and all the variables denote the same physical value with the subscript "l" for the liquid phase and "s" for the solid phase, respectively. The solid phase temperature must be lower than the melting temperature, which serves as one of the conditions for the model validity, as stated in the following. 
\begin{remark}
To keep the physical state of each phase meaningful, the following conditions must be maintained: 
\begin{align}\label{valid1-2ph}
T_{{\rm l}}(x,t) \geq& T_{{\rm m}}, \quad \forall x\in(0,s(t)), \quad \forall t>0, \\
\label{valid2-2ph}T_{{\rm s}}(x,t) \leq& T_{{\rm m}}, \quad \forall x\in(s(t),L), \quad \forall t>0, \\
\label{valid3-2ph} 0< &s(t)<L, \quad \forall t>0. 
\end{align}
\end{remark}
The following assumptions on the initial data \\$(T_{{\rm l},0}(x), T_{{\rm s},0}(x), s_0) := (T_{{\rm l}}(x,0), T_{{\rm s}}(x,0), s(0))$ are imposed.  
\begin{assum}\label{initial-2ph} 
$0<s_0<L$,  $T_{{\rm l},0}(x) \geq T_{{\rm m}}$ for all $x \in [0,s_0]$, $T_{{\rm s},0}(x) \leq T_{{\rm m}}$ for all $x \in [s_0, L]$, and $T_{{\rm l},0}(x) $ and $T_{{\rm s},0}(x) $  are continuously differentiable in $x\in[0,s_0]$ and $x\in[s_0, L]$, respectively. 
 \end{assum}
 \begin{assum}\label{ass:E0}
The initial values satisfy 
\begin{align}
&\frac{k_{{\rm l}}}{\alpha_{{\rm l}}} \int_{0}^{s_0} (T_{{\rm l},0}(x) - T_{{\rm m}}) {\rm d}x  \notag\\
&+ \frac{k_{{\rm s}}}{\alpha_{{\rm s}}} \int_{s_0}^{L} (T_{{\rm s},0}(x) - T_{{\rm m}}) {\rm d}x+ \gm s_0>0 . 
\end{align}
\end{assum}
The following lemma is proven in \cite{Cannon71flux} (page 4, Theorem 1).  
 \begin{lem}\label{solidvalid} 
With Assumption \ref{initial-2ph} and provided that $q_{{\rm c}}(t) \geq 0$ for all $t>0$, there exists $\sigma \in (0, \infty ]$ such that for any $\bar t \leq \sigma $ there is a unique solution of the system \eqref{sys1}--\eqref{sys5} with satisfying the conditions \eqref{valid1-2ph}--\eqref{valid3-2ph} for all $t \in (0, \bar t)$. Moreover, if $\sigma \neq \infty$ then either $s(\sigma) = 0$ or $s(\sigma ) = L$. 
\end{lem}
Furthermore, we impose Assumptions \ref{disturbed} and \ref{ISS-cond}, and the restriction for the setpoint is given as follows. 
\begin{assum}\label{assum2-2ph}
The setpoint is chosen to satisfy
\begin{align} \label{setpoint}
\underline{s}_{{\rm r}} < s_{{\rm r}} < L, 
\end{align} 
where $\underline{s}_{{\rm r}} : = s_0 + \fr{\beta_{{\rm l}}}{\alpha_{{\rm l}}}  \int_{0}^{s_0} (T_{{\rm l},0}(x) - T_{{\rm m}}) {\rm d}x + \fr{\beta_{{\rm s}}}{\alpha_{{\rm s}}} \int_{s_0}^{L} (T_{{\rm s},0}(x) - T_{{\rm m}}) {\rm d}x $ and $\beta_{i} = \frac{k_{i}}{\gamma}$ for $i = {\rm l}, {\rm s}$. 
\end{assum}
\subsection{Control design and an equivalent open-loop analysis} 
We apply the boundary feedback control design developed in \cite{koga2018CDC} 
\begin{align}\label{control-2ph}
q_{ c}(t)=& -c \left( \frac{k_{{\rm l}}}{\alpha_{{\rm l}}} \int_{0}^{s(t)} (T_{{\rm l}}(x,t) - T_{{\rm m}}) {\rm d}x \right. \notag\\
&\left. + \frac{k_{{\rm s}}}{\alpha_{{\rm s}}} \int_{s(t)}^{L} (T_{{\rm s}}(x,t) - T_{{\rm m}}) {\rm d}x + \gamma (s(t) - s_{{\rm r}}) \right), 
\end{align}
to the two-phase Stefan problem with the unknown heat loss $q_{{\rm f}}(t)$ at the boundary of the solid phase governed by \eqref{sys1}--\eqref{sys5}. As in the previous procedure, the equivalence of the closed-loop system under the control law \eqref{control-2ph} with the system under an open-loop input is presented in the following lemma. 
\begin{lem} \label{lem:2ph} 
Under Assumptions \ref{disturbed}, \ref{ISS-cond}, and \ref{initial-2ph}--\ref{assum2-2ph}, the closed-loop system consisting of \eqref{sys1}--\eqref{sys5} with the control law \eqref{control-2ph} has a unique classical solution which is equivalent to the open-loop solution of \eqref{sys1}--\eqref{sys5} with 
\begin{align} \label{open-2ph}
q_{\rm c}(t) = q_{0} e^{-ct} + c \int_0^t e^{-c(t-\tau)} q_{\rm f}(\tau) {\rm d}\tau , 
\end{align}
for all $t \geq 0$, where
\begin{align}\label{q0-2ph}
q_{0} =& - c \left( \frac{k_{{\rm l}}}{\alpha_{{\rm l}}} \int_{0}^{s_0} (T_{{\rm l},0}(x) - T_{{\rm m}}) {\rm d}x \right. \notag\\
&\left. + \frac{k_{{\rm s}}}{\alpha_{{\rm s}}} \int_{s_0}^{L} (T_{{\rm s},0}(x) - T_{{\rm m}}) {\rm d}x + \gamma (s_0 - s_{{\rm r}}) \right). 
\end{align} 
\end{lem} 

\begin{pf} 
Taking the time derivative of the control law \eqref{control-2ph} together with \eqref{sys1}--\eqref{sys5} leads to $\dot{q}_{{\rm c}}(t) = - c(q_{{\rm c}}(t) - q_{{\rm f}}(t))$. The solution to the differential equation is equivalent with the formulation of \eqref{open-2ph} with the initial condition $q_{{\rm c}}(0) = q_0$. The setpoint restriction given in Assumption \ref{assum2-2ph} provides $q_0>0$, and hence $q_{{\rm c}}(t)>0$ for all $t \geq 0$. Applying Lemma \ref{solidvalid}, we ensure the existence of $\sigma \in (0, \infty]$ stated in Lemma \ref{solidvalid}. Moreover, we prove $\sigma = \infty$ by contradiction, similarly to the proof of Lemma \ref{global}. Assume the existence of $\sigma>0$ such that the conditions \eqref{valid1-2ph}--\eqref{valid3-2ph} are satisfied for $ t \in (0, \sigma)$. Taking into account the specific heat in the solid phase, the internal energy of the total system is defined by 
\begin{align}\label{Et-2ph}
E(t) =& \frac{k_{{\rm l}}}{\alpha_{{\rm l}}} \int_{0}^{s(t)} (T_{{\rm l}}(x,t) - T_{{\rm m}}) {\rm d}x \notag\\
&+ \frac{k_{{\rm s}}}{\alpha_{{\rm s}}} \int_{s(t)}^{L} (T_{{\rm s}}(x,t) - T_{{\rm m}}) {\rm d}x + \gm s(t). 
\end{align} 
Taking the time derivative of \eqref{Et-2ph} leads to the energy conservation law given by the same formulation as \eqref{energy-1ph} for the case of one-phase Stefan problem. In the similar manner as the proof of Lemma \ref{global}, we deduce 
\begin{align} \label{Epositive} 
E(t)>0, \quad \forall t \geq 0. 
\end{align} 
 If $s(\sigma) = 0$, then substituting it into \eqref{Et-2ph} with the help of \eqref{valid2-2ph} leads to $E(\sigma) <0$, which contradicts with \eqref{Epositive}. Furthermore, applying $q_{{\rm c}}(t) >0$ and the condition \eqref{valid1-2ph} for $t \in (0, \sigma)$ to the control law \eqref{control-2ph}, we get 
 \begin{align} \label{Enegative} 
 \frac{k_{{\rm s}}}{\alpha_{{\rm s}}} \int_{s(t)}^{L} (T_{{\rm s}}(x,t) - T_{{\rm m}}) {\rm d}x + \gamma (s(t) - s_{{\rm r}})<0,
 \end{align} 
for all $t \in (0, \sigma)$, which ensures at least $s(\sigma) < L$. Hence, we conclude that $\sigma = \infty$, which guarantees the well-posedness of the closed-loop solution and the conditions for the model validity \eqref{valid1-2ph}--\eqref{valid3-2ph} to hold globally, i.e., for all $t \geq 0$, which completes the proof of Lemma \ref{lem:2ph}.  

\end{pf} 

\subsection{Main result and proof} 
 We present ISS result for the closed-loop system in the following theorem. 
\begin{thm}\label{theo-2}
Under Assumptions \ref{disturbed}, \ref{ISS-cond}, and \ref{initial-2ph}--\ref{assum2-2ph}, the closed-loop system consisting of \eqref{sys1}--\eqref{sys5} with the control law \eqref{control-2ph} satisfies the model validity conditions \eqref{valid1-2ph}--\eqref{valid3-2ph}, and is ISS with respect to the heat loss $q_{{\rm f}}(t)$ at the boundary of the solid phase. Moreover, there exist positive constants $M_1>0$ and $M_2>0$ such that the following estimate holds:
\begin{align}\label{h1-2ph}
\Psi(t) \leq M_1 \Psi(0) e^{- \lambda t} + M_2  \sup_{\tau \in [0,t]} |q_{{\rm f}}(\tau)| ,
\end{align}
for all $t\geq0$, where $\lambda = \frac{1}{32} \min \left\{ \fr{\alpha_{{\rm l}}}{ L^2}, \fr{2 \alpha_{{\rm s}}}{ L^2}, c \right\} $, in the $L_2$ norm \\
$\Psi(t) = \left( \int_{0}^{s(t)} (T_{{\rm l}}(x,t) - T_{{\rm m}})^2 {\rm d}x \right. \\
\left. + \int_{s(t)}^{L} (T_{{\rm s}}(x,t) - T_{{\rm m}})^2 {\rm d}x + (s(t) - s_{{\rm r}})^2 \right)^{\frac{1}{2}}$.  
\end{thm}
\begin{pf} 
The proof of the conditions \eqref{valid1-2ph}--\eqref{valid3-2ph} are given in the proof of Lemma \ref{lem:2ph}. To prove ISS, by following the procedure in \cite{koga2018CDC}, first we introduce the reference error states as follows. 
\begin{align}\label{ref-2ph}
u(x,t) :=& T_l (x,t)-T_{{\rm m}}, \\
X(t) :=& s(t) - s_{{\rm r}} + \fr{\beta_{{\rm s}}}{\alpha_{{\rm s}}} \int_{s(t)}^{L} (T_{{\rm s}}(x,t) - T_{{\rm m}}) {\rm d}x . \label{refX-2ph} 
\end{align}
Using these reference error variables, the total PDE-ODE-PDE system given in \eqref{sys1}--\eqref{sys5} is reduced to the following PDE-ODE system
\begin{align}\label{uX-sys1}
 u_{t}(x,t) =& \alpha_{{\rm l}} u_{xx}(x,t), \quad 0<x<s(t), \\
\label{uX-BC}u_x(0,t) = &- q_{{\rm c}}(t)/k_{{\rm l}}, \quad  u(s(t),t) = 0,\\
\dot{X}(t) =& - \beta_{{\rm l}} u_x(s(t),t) - d(t),  \label{uX-sys3}
\end{align}
where $d(t) = \fr{\beta}{k} q_{\rm f}(t)$. Note that the formulation of the above system is equivalent to \eqref{u-sys1}--\eqref{u-sys4} which is the reference error system in the one-phase case.  Therefore, following the same procedure as in Section \ref{sec:closed}, it is straightforward to derive that there exist positive constants $N_1>0$ and $N_2>0$ such that the following norm estimate holds: 
\begin{align} \label{phiineq-2ph}
\Phi(t) \leq & N_1 \Phi(0) e^{- \frac{b}{4}t} + N_2 \sup_{\tau \in [0, t]} | q_{{\rm f}}(\tau) |, 
\end{align} 
where $\Phi(t) = \left( \int_{0}^{s(t)} u(x,t)^2 {\rm d}x + X(t)^2 \right)^{\frac{1}{2}}$, $ b =\fr{1}{8} \min \left\{ \fr{\alpha_{{\rm l}}}{L^2}, c \right\} $. 
Let us define the following three functionals 
\begin{align}\label{V1def-2ph}
V_1 (t)=& \int_0^{s(t)}  (T_{{\rm l}}(x,t) - T_{{\rm m}})^2 {\rm d}x, \\
V_2 (t)=& \int_{s(t)}^{L} (T_{{\rm s}}(x,t) - T_{{\rm m}})^2 {\rm d}x, \label{V2def-2ph} \\
V_3 (t)=& (s(t) - s_{{\rm r}})^2. \label{V3def-2ph}
\end{align} 
Taking the time derivative of \eqref{V2def-2ph} along with the solid phase dynamics \eqref{sys3} and \eqref{sys4}, we get
\begin{align} 
 \dot{V}_2 (t)= & -2 \alpha_s \int_{s(t)}^{L} \left(\frac{\pa T_{{\rm s}}}{\pa x}(x,t) \right)^2 {\rm d}x \notag\\
&- \frac{2 \alpha_s}{k_s} (T_{{\rm s}}(L,t) - T_{{\rm m}}) q_{{\rm f}}(t) . \label{V2dot} 
\end{align} 
Applying Young's, Cauchy-Schwarz, Poincare's and Agmon's inequalities to \eqref{V2dot}, we arrive at the following differential inequality 
\begin{align} \label{Vdot2ineq-2ph} 
\dot{V}_2 (t)\leq - \frac{\alpha_{{\rm s}}}{4 L^2} V_2(t) + \frac{4 L \alpha_{{\rm s}}}{k_{{\rm s}}^2} q_{{\rm f}}(t)^2  . 
\end{align} 
Applying the same procedure as the derivation from \eqref{Wdot} to \eqref{Wineq}, and taking the square root, one can derive 
\begin{align} \label{V2ineq-2ph}
\sqrt{V_2 (t)} \leq \sqrt{V_2(0)} e^{-  \frac{\alpha_{{\rm s}}}{16 L^2} t} + \frac{4L \sqrt{2L}}{k_{{\rm s}}} \sup_{\tau \in [0, t]} |q_{{\rm f}}(\tau) |  . 
\end{align} 
Taking the square of \eqref{refX-2ph}, and applying Young's and Cauchy-Schwarz inequalities with the help of $0 < s(t) < L$, one can obtain the following inequality, 
\begin{align} \label{Xineq}
 X(t)^2 \leq & 2 V_3(t) + \fr{2 L \beta_{{\rm s}}^2}{\alpha_{{\rm s}}^2} V_2 (t). 
\end{align}
Applying the same manner to the relation $s(t) - s_{{\rm r}} = X(t) - \frac{\beta_s}{\alpha_s} \int_{s(t)}^{L} (T_{{\rm s}}(x,t) - T_{{\rm m}}) {\rm d}x$ obtained by \eqref{refX-2ph}, one can also derive 
\begin{align}\label{Ytineq}
 V_3(t)   \leq  2 X(t)^2 + \fr{2 L \beta_{{\rm s}}^2}{\alpha_{{\rm s}}^2} V_2 (t). 
\end{align}
Combining \eqref{phiineq-2ph}, \eqref{V2ineq-2ph}, \eqref{Xineq}, and \eqref{Ytineq} using the definitions in \eqref{V1def-2ph}--\eqref{V3def-2ph}, the estimate of the norm $\Psi(t) = \sqrt{V_1(t) + V_2(t) + V_3(t)}$ is obtained by the inequality \eqref{h1-2ph} for some positive constants $M_1>0$ and $M_2>0$, which completes the proof of Theorem \ref{theo-2}. 
\end{pf} 

\section{Numerical Simulation}\label{simulation}
	\begin{table}[t]
	\vspace{2mm}
	
\caption{Physical properties of zinc}
\begin{center}
    \begin{tabular}{| l | l | l | }
    \hline
    $\textbf{Description}$ & $\textbf{Symbol}$ & $\textbf{Value}$ \\ \hline
    Density & $\rho$ & 6570 ${\rm kg}\cdot {\rm m}^{-3}$\\ 
    Latent heat of fusion & $\Delta H^*$ & 111,961${\rm J}\cdot {\rm kg}^{-1}$ \\ 
    Heat Capacity & $C_p$ & 389.5687 ${\rm J} \cdot {\rm kg}^{-1}\cdot {\rm K}^{-1}$  \\  
    Thermal conductivity & $k$ & 116 ${\rm w}\cdot {\rm m}^{-1}$  \\ \hline
    \end{tabular}
\end{center}
\end{table}
\begin{figure}[t]
\centering 
\subfloat[Convergence of the interface with respect to the maximum value of heat flux without the overshoot.]
{\includegraphics[width=2.7in]{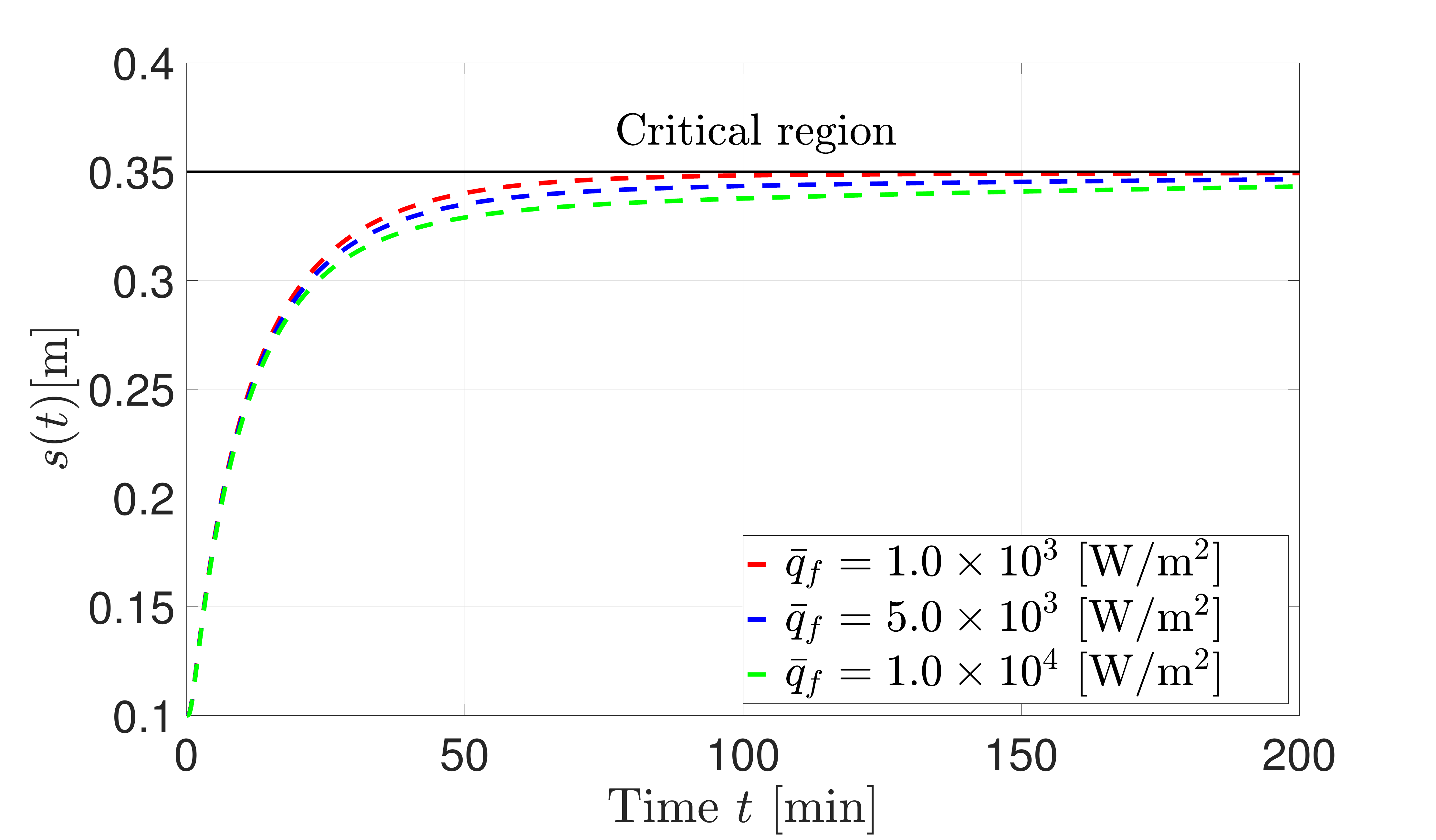}\label{fig:interface}}\\
\subfloat[Positivity of the closed-loop controller maintains.]
{\includegraphics[width=2.7in]{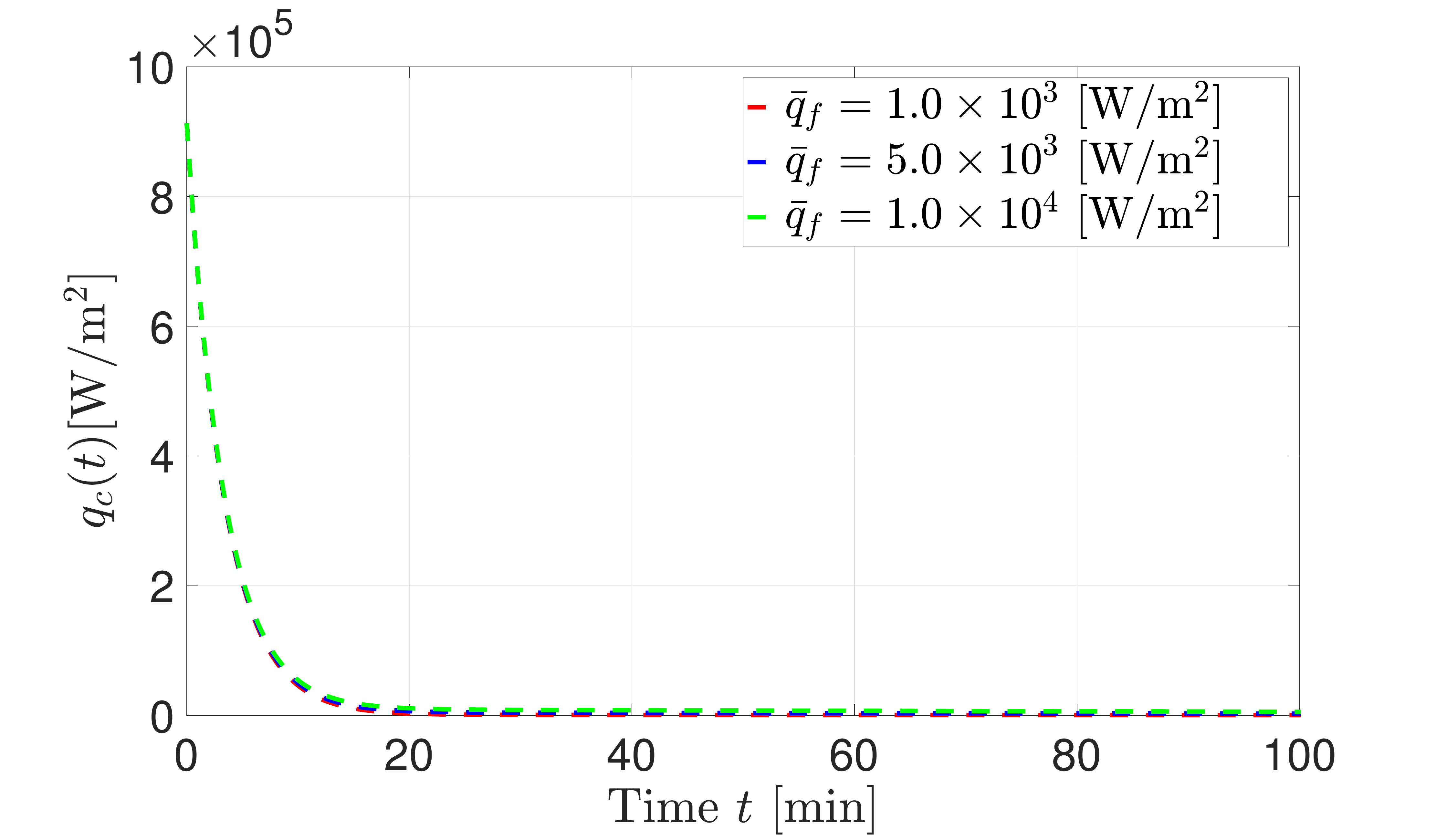}\label{fig:qc}}\\
\subfloat[The model validity of the boundary liquid temperature holds, i.e., $T(0,t)>T_{{\rm m}}$.]
{\includegraphics[width=2.7in]{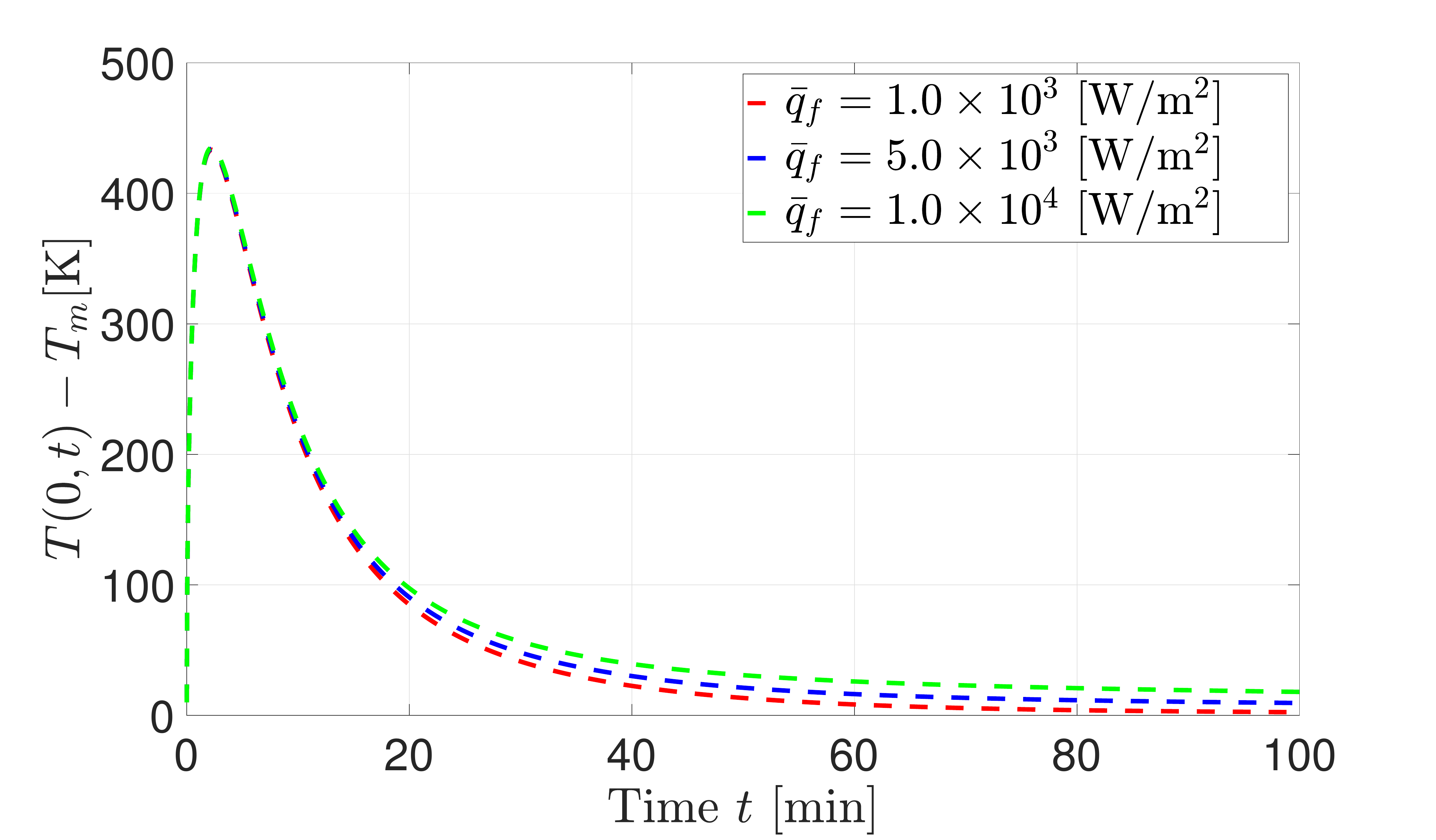}\label{fig:T0}}
\caption{The responses of the system \eqref{eq:stefanPDE}-\eqref{eq:stefanODE} with \eqref{control}.  }
\label{response}
\end{figure}

Simulation results are performed for the one-phase Stefan problem with the heat loss at the interface, by considering a strip of zinc as in \cite{maidi2014,Shumon16} whose physical properties are given in Table 1. Here, we use the well known boundary immobilization method combined with finite difference semi-discretization \cite{kutluay97}. The setpoint and the initial values are chosen as $s_{{\mathrm r}}$ = 0.35 m, $s_0$ = 0.1 m, and $T_0(x)-T_{{\mathrm m}}= \bar{T}_0(1-x/s_0)$ with $ \bar{T}_0$ = 10 K. Then, the setpoint restriction \eqref{eq:setpoint} is satisfied. The control gain is set as $c =$ 5.0 $\times$ 10$^{-3}$/s, and the heat loss at the interface is set as 
\begin{align} 
q_{\rm f}(t) = \bar{q}_{\rm f} e^{- Kt}, 
\end{align} 
where $K$ = 5.0 $\times$ 10$^{-6}$ /s.  The closed-loop responses for $\bar{q}_{\rm f}$ = 1.0 $\times $ 10$^3$W/m$^2$ (red), 5.0 $\times $ 10$^3$W/m$^2$ (blue), and 10.0 $\times $ 10$^3$W/m$^2$ (green) are implemented as depicted in Fig \ref{fig:interface}-\ref{fig:T0}. Fig \ref{fig:interface} shows the dynamics of the interface, which illustrates the convergence to the setpoint with an error due to the unknown heat loss at the interface. This error becomes larger as $\bar{q}_{{\rm f}}$ gets larger, which is consistent with the ISS result. In addition, the property $0<s(t)<s_{{\rm r}}$ is observed, which is consistent with Lemma \ref{closed-validity}. Fig \ref{fig:qc} shows the dynamics of the proposed closed-loop control law, and Fig \ref{fig:T0} shows the dynamics of the boundary temperature $T(0,t)$. These figures illustrate the other conditions proved in Lemma \ref{closed-validity}. Hence, we can observe that the simulation results are consistent with the theoretical result we prove as model validity conditions and the ISS sense.

\section{Conclusions}\label{conclusion}
In this paper we study the Stefan problem with an unknown heat loss. First we revisit the one-phase Stefan problem with heat loss at the interface studied in \cite{koga2018ISS}. The closed-loop control law is designed based on \cite{Shumon16}, and the well-posedness and the conditions for the model validity are proved based on the open-loop analysis. A nonlinear backstepping transformation for moving boundary problems is utilized, and the estimate of the $L_2$ norm is derived in the sense of ISS with respect to the heat loss at the interface. Furthermore, the extension of the ISS result to the two-phase Stefan problem with the heat loss at the boundary of the solid phase under the control law proposed in \cite{koga2018CDC} is achieved. 


\bibliographystyle{unsrt}

\begin{thebibliography}{9} 

\bibitem{arcak01}
M. Arcak and P. Kokotovic, 
\newblock ``Nonlinear observers: a circle criterion design and robustness analysis,"
\newblock {\em Automatica}, vol. 37, pp.1923-1930, 2001.

\bibitem{Armaou01}
A. Armaou and P.D. Christofides,
\newblock ``Robust control of parabolic PDE systems with time-dependent spatial domains,''
\newblock {\em Automatica}, vol. 37, pp. 61--69,
  2001.
  
  \bibitem{Leobardo2018}
  L. Camacho-Solorio, S. Moura, and M. Krstic, 
  \newblock ``Robustness of boundary observers for radial diffusion equations to parameter uncertainty,"
  \newblock {\em 2018 American Control Conference (ACC)}, pages 3484-3489.
  IEEE, 2018.
  
    \bibitem{Cannon71flux}
J. R. Cannon and M. Primicerio, 
\newblock ``A two phase Stefan problem with flux boundary conditions," 
\newblock {\em Annali di Matematica Pura ed Applicata,} 88.1, 193-205, 1971.

\bibitem{Chen2008} 
X. Chen, J. Chadam, L. Jiang, and W. Zheng, 
\newblock ``Convexity of the exercise boundary of the American put option on a zero dividend asset," 
\newblock {\em Mathematical Finance: An International Journal of Mathematics, Statistics and Financial Economics}, vol. 18, no. 1, pp.185-197, 2008 
  
  \bibitem{chung2004} 
H. Chung, and S. Dans, 
\newblock ``Numerical modeling of scanning laser-induced melting, vaporization and resolidification in metals subjected to step heat flux input,"
\newblock {\em International journal of heat and mass transfer}, vol. 47, pp. 4153-4164, 2004. 



\bibitem{conrad_90}
F. Conrad, D. Hilhorst, and T.I. Seidman,
\newblock ``Well-posedness of a moving boundary problem arising in a
 dissolution-growth process,''
\newblock {\em Nonlinear Analysis}, vol. 15, pp. 445--465, 1990.

\bibitem{Daraoui2010}
N.~Daraoui, P.~Dufour, H.~Hammouri, and A.~Hottot,
\newblock ``Model predictive control during the primary drying stage of
  lyophilisation,"
\newblock {\em Control Engineering Practice}, vol. 18, pp. 483--494, 2010.

\bibitem{Du2010speading} 
Y. Du, and Z. Lin, 
\newblock ``Spreading-vanishing dichotomy in the diffusive logistic model with a free boundary,"
\newblock {\em SIAM Journal on Mathematical Analysis}, vol. 42(1), pp. 377--405, 2010.





\bibitem{Friedman59}
A. Friedman,
\newblock ``Free boundary problems for parabolic equations I. Melting of solids," 
\newblock  {\em Journal of Mathematics and Mechanics}, 8(4), pp.499-517, 1959.

\bibitem{Friedman1999}
A. Friedman and F. Reitich, 
\newblock ``Analysis of a mathematical model for the growth of tumors,"
\newblock  {\em Journal of mathematical biology}, vol. 38(3), pp.262-284, 1999.

\bibitem{Gupta03}
S.~Gupta,
\newblock {\em The Classical Stefan Problem. Basic Concepts, Modelling and
  Analysis}.
\newblock North-Holland: Applied mathematics and Mechanics, 2003.


\bibitem{Karafyllis16}
I. Karafyllis, and M. Krstic, 
\newblock ``ISS with respect to boundary disturbances for 1-D parabolic PDEs"
\newblock {\em IEEE Transactions on Automatic Control}, 61(12), pp.3712-3724, 2016. 

\bibitem{Karafyllis17}
I. Karafyllis, and M. Krstic, 
\newblock ``ISS in different norms for 1-D parabolic PDEs with boundary disturbances"
\newblock {\em SIAM Journal on Control and Optimization}, 55(3), pp.1716-1751, 2017. 

  \bibitem{Karafyllis-issbook}
I. Karafyllis, and M. Krstic, 
\newblock {\em Input-to-state stability for PDEs}.
\newblock Springer-Verlag, London (Series: Communications and Control Engineering), 2019.


\bibitem{Shumon16}
S.~Koga, M.~Diagne, S.~Tang, and M.~Krstic,
\newblock ``Backstepping control of the one-phase stefan problem,"
\newblock In {\em 2016 American Control Conference (ACC)}, pages 2548--2553.
  IEEE, 2016.
  
  \bibitem{Shumon16CDC}
S. Koga, M. Diagne, and M. Krstic,
\newblock ``Output feedback control of the one-phase Stefan problem,"
\newblock In {\em 55th Conference on Decision and Control (CDC)}, pages 526--531.
  IEEE, 2016.

  
  \bibitem{Shumon19journal}
 S. Koga, M. Diagne, and M. Krstic,
  \newblock ``Control and state estimation of the one-phase Stefan problem via backstepping design,"
  \newblock {\em IEEE Transactions on Automatic Control}, vol. 64, no. 2, pp. 510--525, 2019. 
  
    
  \bibitem{Shumon2017ACC}
S. Koga, R. Vazquez and M. Krstic, 
\newblock ``Backstepping control of the Stefan problem with flowing
liquid,''
\newblock In {\em 2017 American Control Conference (ACC)}, pages
1151-1156. IEEE 2017

\bibitem{Shumon17seaice}
S. Koga and M. Krstic, 
\newblock ``Arctic sea ice temperature profile estimation via backstepping observer design," 
\newblock In {\em 2017 Conference on Control Technology and Applications (CCTA)}, pages 1722-1727, IEEE, 2017.

  \bibitem{koga2017battery}
S. Koga, L. Camacho-Solorio, and M. Krstic, 
\newblock ``State Estimation for Lithium Ion Batteries With Phase Transition Materials," 
\newblock |n {\em ASME 2017 Dynamic Systems and Control Conference}, American Society of Mechanical Engineers, 2017.

  \bibitem{koga2017CDC} 
S. Koga and M. Krstic, 
\newblock ``Delay compensated control of the Stefan problem," 
\newblock In {\em 56th Conference on Decision and Control (CDC)}, pp. 1242-1247, IEEE, 2017.


  \bibitem{koga2018polymer} 
S. Koga, D. Straub, M. Diagne, and M. Krstic,
\newblock ``Thermodynamic Modeling and Control of Screw Extruder for 3D Printing,"
\newblock In {\em 2018 American Control Conference (ACC)}, pages 2551--2556.
  IEEE, 2018.
  
  \bibitem{koga2018ISS} 
S. Koga, I. Karafyllis, and M. Krstic,
\newblock ``Input-to-State Stability for the Control of Stefan Problem with Respect to Heat Loss at the Interface,"
\newblock In {\em 2018 American Control Conference (ACC)}, pages 1740--1745.
  IEEE, 2018.

  
      \bibitem{koga2018CDC} 
S. Koga and M. Krstic, 
\newblock ``Control of Two-Phase Stefan Problem via Single Boundary Heat Input," 
\newblock In {\em 57th Conference on Decision and Control (CDC)}, pp. 2914-2919, IEEE, 2018.

\bibitem{koga_2019delay}
S. Koga, D. Bresch-Pietri, and M. Krstic, 
\newblock ``Delay compensated control of the Stefan problem and robustness to delay mismatch," 
\newblock Preprint, available at \\ http://arxiv.org/abs/1901.09809, 2019. 


  
  \bibitem{krstic2008boundary}
M.~Krstic and A.~Smyshlyaev,
\newblock {\em Boundary Control of {PDE}s: A Course on Backstepping Designs}.
\newblock Singapore: SIAM, 2008.

\bibitem{krstic2009delay}
M. Krstic, 
\newblock {\em Delay compensation for nonlinear, adaptive, and PDE systems}.
\newblock Birkh{\"a}user Boston, 2009.


\bibitem{kutluay97}
S. Kutluay, A. R. Bahadir, and A. {\"O}zdes, 
\newblock ``The numerical solution of one-phase classical Stefan problem," 
\newblock {\em Journal of computational and applied mathematics}, 81.1, pp. 135-144, 1997.

\bibitem{Lei2013}
C. Lei, Z. Lin, and H. Wang, 
\newblock ``The free boundary problem describing information diffusion in online social networks,"
\newblock {\em Journal of Differential Equations}, vol. 254(3), pp.1326-1341, 2013.


\bibitem{maidi2014}
A.~Maidi and J.-P. Corriou,
\newblock ``Boundary geometric control of a linear stefan problem,"
\newblock {\em Journal of Process Control}, vol. 24, pp. 939--946, 2014.

\bibitem{mcgilly2015}
L.J. McGilly, P. Yudin, L. Feigl, A.K. Tagantsev, and N. Setter,
\newblock ``Controlling domain wall motion in ferroelectric thin films,"
\newblock {\em Nature nanotechnology}, vol. 10(2), pp. 145, 2015.

\bibitem{Moura2014} 
S.J. Moura, N.A. Chaturvedi, and M. Krstic, 
\newblock ``Adaptive partial differential equation observer for battery state-of-charge/state-of-health estimation via an electrochemical model," 
\newblock {\em Journal of Dynamic Systems, Measurement, and Control}, vol. 136, no. 1, pp. 011015-1?011015-11, 2014. 

  \bibitem{petrus2012}
B. Petrus, J. Bentsman, and B.G. Thomas,
\newblock ``Enthalpy-based feedback control algorithms for the Stefan problem,"
\newblock {\em Decision and Control (CDC), 2012 IEEE 51st Annual Conference on}, pp. 7037--7042, 2012.

\bibitem{Rabin1998}
Y. Rabin, and A. Shitzer,
\newblock ``Numerical solution of the multidimensional freezing problem during cryosurgery,"
\newblock {\em Journal of biomechanical engineering}, 120(1), 32-37, 1998.



\bibitem{Sherman67}
B. Sherman, 
\newblock ``A free boundary problem for the heat equation with prescribed flux at both fixed face and melting interface", 
\newblock {\em Quarterly of Applied Mathematics}, 25(1), pp.53-63, 1967.


\bibitem{andrew2004}
A. Smyshlyaev and M. Krstic,
\newblock ``Closed-form boundary State feedbacks for a class of 1-D partial integro-differential equations,"
\newblock {\em IEEE Transactions on Automatic Control,}, vol. 49, pp. 2185--2202, 2004.

\bibitem{Sontag08} 
E.D. Sontag, 
\newblock ``Input to state stability: Basic concepts and results,"
\newblock In Nonlinear and optimal control theory (pp. 163-220). Springer Berlin Heidelberg, 2008.



\bibitem{zalba03}
B. Zalba, J.M. Marin, L.F. Cabeza, and H. Mehling,
\newblock ``Review on thermal energy storage with phase change: materials, heat transfer analysis and applications,"
\newblock {\em Applied thermal engineering}, vol. 23, pp. 251--283, 2003.



\end{thebibliography}

\end{document}